\documentclass[12pt,reqno]{amsart}

\usepackage{hyperref,setspace}
\usepackage{fullpage}

\numberwithin{equation}{section}
\def\e{\epsilon}

\def\bp{\begin{proposition}}
\def\ep{\end{proposition}}
\def\bt{\begin{theo}}
\def\et{\end{theo}}
\def\be{\begin{equation}}
\def\ee{\end{equation}}
\def\bl{\begin{lemma}}
\def\el{\end{lemma}}
\def\bc{\begin{corollary}}
\def\ec{\end{corollary}}
\def\pr{\noindent{\bf Proof: }}

\def\bd{\begin{definition}}
\def\ed{\end{definition}}

\def\max{{\rm max\,}}

\def\max{{\rm max\,}}

\newtheorem{theo}{Theorem}[section]
\newtheorem{lemma}{Lemma}[section]
\newtheorem{definition}{Definition}[section]
\newtheorem{corollary}{Corollary}[section]
\newtheorem{proposition}{Proposition}[section]
\theoremstyle{theorem}
\newtheorem{thm}{Theorem}[section]

\theoremstyle{definition}

\theoremstyle{remark}
\newtheorem{rem}[thm]{Remark}

\newcommand{\N}{\mathbb N}

\begin{document}

\title{Sampling, Metric Entropy and Dimensionality Reduction}

\author{Dima Batenkov}
\address{Department of Mathematics, The Weizmann Institute of Science, Rehovot 76100, Israel.}
\email{dima.batenkov@gmail.com}
\thanks{}

\author{Omer Friedland}
\address{Institut de Math\'ematiques de Jussieu, Universit\'e Pierre et Marie Curie (Paris 6), 4 Place Jussieu, 75005 Paris, France.}
\email{friedland@math.jussieu.fr}
\thanks{}

\author{Yosef Yomdin}
\address{Department of Mathematics, The Weizmann Institute of Science, Rehovot 76100, Israel.}
\email{yosef.yomdin@weizmann.ac.il}
\thanks{This research was supported by the Adams Fellowship Program of the Israel Academy of Sciences and Humanities, by the ISF, Grant No. 639/09 and by the Minerva foundation.}

\begin{abstract}
Let $Q$ be a relatively compact subset in a Hilbert space $V$. For a given $\e>0$ let $N(\e,Q)$ be the minimal number of linear measurements, sufficient to reconstruct any $x \in Q$ with the accuracy $\e$. We call $N(\e,Q)$ a sampling $\e$-entropy of $Q$. Using Dimensionality Reduction, as provided by the Johnson-Lindenstrauss lemma, we show that, in an appropriate probabilistic setting, $N(\e,Q)$ is bounded from above by the Kolmogorov's $\e$-entropy $H(\e,Q)$, defined as $H(\e,Q)=\log M(\e,Q)$, with $M(\e,Q)$ being the minimal number of $\e$-balls covering $Q$. As the main application, we show that piecewise smooth (piecewise analytic) functions in one and several variables can be sampled with essentially the same accuracy rate as their regular counterparts. For univariate piecewise $C^k$-smooth functions this result, which settles the so-called Eckhoff conjecture, was recently established in \cite{Bat} via a deterministic ``algebraic reconstruction'' algorithm.
\end{abstract}

\maketitle

\section{Introduction}

In this paper we study the ``sampling $\e$-entropy'' of functional classes in a Hilbert space $V$.
``Samples'' or ``measurements'' of an element $x\in V$ are assumed to be linear functionals on $V$, i.e. the values of
the scalar products
\begin{align}\label{scalar}
\lambda_v (x)=\langle x,v \rangle,
\end{align}
where $v\in V$ is a known vector. To simplify the presentation we shall assume in this
paper that the measurements are exact. See Remark \ref{fin.acc}, where we discuss measurements of a finite accuracy.

\smallskip

Let $Q$ be a relatively compact subset in $V$. We start with Kolmogorov's metric entropy, which measures ``information complexity''
of $Q$ in a prescribed resolution $\e$.

\bd\label{Kolm} For each $\e>0$ the covering number $M(\e,Q)$ is the minimal number of $\e$-balls $B_\e$ in $V$ covering $Q$.
$H(\e,Q)=\log_2 M(\e,Q)$ is called the Kolmogorov $\e$-entropy of $Q$.
\ed
Metric entropy $H(\e,Q)$ of functional classes $Q$ was introduced and studied in \cite{Kol1,Kol2,Kol.Tih}. It is exactly the number
of binary bits of information we need to describe any element $x\in Q$ with the accuracy $\e$. Indeed, to get such a description it
is necessary and sufficient to specify the ball $B_\e$ (or its center) in the minimal $\e$-covering of $Q$, containing $x$. The number
of the balls is $M(\e,Q)$, and hence the number of bits required is $H(\e,Q)=\log_2 M(\e,Q)$.

\smallskip

In our definition of the ``sampling $\e$-entropy'' $N(\e,Q)$ we just add the requirement that the information used to reconstruct
unknown $x\in Q$ be provided by linear measurements $\lambda_v(x)$ as defined in \eqref{scalar}.

\bd\label{Sampl.Entr}
For each given $\e >0$ the sampling $\e$-entropy $N(\e,Q)$ of $Q$ is defined as the minimal number of exact
linear measurements $\lambda_v(x)$, sufficient to reconstruct any $x \in Q$ with the accuracy $\e$. The set of measurements
used may depend on $Q$ and on $\e$, but not on the specific $x\in Q$ to be reconstructed.
\ed
Below we recall some classical results providing sampling entropy for $C^k$-smooth and analytic functions. Based on a recent result of
\cite{Bat} we obtain a sharp bound for sampling entropy for univariate piecewise $C^k$-smooth functions. However, mostly
in this paper we shall work with a randomized version of sampling entropy. Assume that for each $N$ a procedure $\mathcal P$ of a random
choice of $N$ measurements (i.e. of an orthonormal system of $N$ vectors in $V$) has been fixed. Let us stress that this procedure depends
only on the ambient Hilbert space $V$ and on the number of measurements $N$. One specific choice of such a procedure $\mathcal P$ is described
in Section \ref{JLL}, where we present a version of Dimensionality Reduction, suitable for our purposes, as provided by the
Johnson-Lindenstrauss lemma (\cite{Joh.Lin}).
\bd\label{Ran.Sampl.Entr}
Let $Q$ be a relatively compact subset in a Hilbert space $V$. For each given $0<p<1$ and for each $\e >0$
the sampling $(\e,p)$-entropy $N(\e,p,Q)$ of $Q$ is defined as the minimal number $N$ with the following property: For each fixed $x\in Q$
apply the procedure $\mathcal P$ to pick randomly $N$ linear measurements $\lambda_{v_j}, \ j=1,\dots,N$. Then with probability at
least $p$ the element $x \in Q$ can be reconstructed from the measurements values $\lambda_{v_j}(x), \ j=1,\dots,N$, with the accuracy $\e$.
\ed
Notice that the random choice according to the procedure $\mathcal P$ in Definition \ref{Ran.Sampl.Entr} does not depend on a specific $x\in Q$. 
Hence it can be performed ``off line'', once forever, and
still we can expect that with probability $p$ any given $x\in Q$ can be $\e$-reconstructed from the chosen set of measurements. However, this last
assertion can be now interpreted in a different way, referring to a choice of $x$, as the measurements have been fixed in advance. In particular,
one can expect that it is also true that for a randomly chosen $x\in Q$ (with respect to a certain natural probability distribution) with
probability $p$ this $x$ can be $\e$-reconstructed from the chosen set of measurements.

\smallskip

The Johnson-Lindenstrauss lemma is a finite-dimensional result, so in order to apply it we have to find for each $\e>0$ a sufficiently
accurate (on $Q$) finite-dimensional approximation of the Hilbert space $V$. We achieve it by making the following assumption, which is satisfied in
all the specific examples considered below:

\smallskip

\noindent{\bf Assumption A} {\it There exists an orthonormal basis $(v_1,v_2,\ldots)$ in $V$ such that for each $d\in\N$ the subspace $V_d \subset V$ 
spanned by $v_1,v_2,\ldots,v_d$ satisfies the following: For any $x\in Q$ we have
$$
|x-\pi_{V_d}(x)|_2\le Cd^{-\beta}|x|_2 ,
$$
where $\beta> 0$, $C>0$, $\pi_{V_d}$ is the projection of $V$ on $V_d$, and $|\cdot|_2$ is the Euclidean norm on $V$.}

\smallskip

For example, Assumption A holds for $Q$ being the class of univariate piecewise-$C^k$ functions, and the subspaces $V_d$ of $d$-th partial sums of 
the Fourier basis of $L^2([-\pi,\pi])$, with $\beta = 1$.

\smallskip

Our main result is the following:

\bt \label{main.short}
Let $Q$ be a relatively compact subset in a Hilbert space $V$ satisfying Assumption A. Then for each $\e>0$ and $0<p<1$,

$$
N(\e,p,Q) \le {{20}\over {1-p}} H\bigl({\e \over 6},Q\bigr).
$$
\et
The proof is given in Section \ref{Entr.JLL}.

\smallskip

Notice that the parameters $C,\beta$ in the Assumption A for $Q,V$ {\it do not enter} the statement of Theorem \ref{main.short}. However, these parameters affect the required measurements accuracy (compare to Remark \ref{fin.acc}).

So, at least for exact random measurements, the requirement that the information on $x\in Q$ comes from {\it linear samples} does not increase significantly the number of samples we need in order to achieve the prescribed reconstruction accuracy.

\smallskip

As the main application, we show in Section \ref{EXa} that {\it multivariate and univariate piecewise smooth functions, as well
as univariate piecewise analytic ones, can be sampled with essentially the same accuracy rate as their regular counterparts}.
This fact was conjectured in various forms in many publications in the field. See discussion after Theorem \ref{pCke}.

\smallskip

As a direct consequence of the main result of this paper, we confirm a certain very general principle, which compares sampling
entropy of the classical functional classes, defined by certain regularity assumptions, and of their extensions, which preserve complexity,
but may completely destroy regularity. Examples considered in Section \ref{EXa} include piecewise-regular functions, finite-dimensional
``warping'', and adding finite-dimensional non-regular signals. In all these cases {\it the randomized sampling entropy of the extended
non-regular classes remains essentially the same, as for their regular counterparts}.

\subsection{Some references and acknowledgement}\label{ack}

Several results on Compressed Sensing and its modifications have recently appeared, where the Johnson-Lindenstrauss lemma,
as well as Kolmogorov's $n$-width, play an important role. See
\cite{Bar.Dav.DeV.Wak,Bar.Cev.Dua.Heg,Can,Coh.Dah.DeV,Dua.Eld,Kas.Tem} and references therein. On the other hand, recent results in
\cite{Ver1} provide important connections of the sampling entropy of subsets with their ``mean width'' and other
geometric characteristics. It would be interesting to compare these results with those of the present paper.

\smallskip

The authors would like to thank A. Melnikov and G. Schechtman for useful discussions, which, in particular, allowed us to correct
some steps in the proof, and to significantly improve the presentation.

\section{Johnson-Lindenstrauss lemma} \label{JLL}

This result (\cite{Joh.Lin}; see also \cite{Bar.Dav.DeV.Wak,Bar.Cev.Dua.Heg,Can,Coh.Dah.DeV,Dua.Eld,Kas.Tem} and references therein)
is central in numerous problems involving dimensionality reduction. We need the following special case of it:

\bt\label{JLLT}
Let $Z\subset {\mathbb R}^d$ be a finite set of cardinality $m$. Put $n={20\over {1-p}} \ln m$, with $p$ being a given probability in
$(0,1)$. Let $W$ be a randomly chosen linear $n$-dimensional subspace in ${\mathbb R}^d$ and let $\pi_W: {\mathbb R}^d \rightarrow W$
be the orthogonal projection onto $W$. Then with probability at least $p$ for each $x,y\in Z$ we have

\be\label{jll}
{1\over 2}|x-y|_2 \le \sqrt {d\over n} \ |\pi_W(x) - \pi_W(y)|_2 \le 2|x-y|_2.
\ee
\et
We recall (following \cite{Ail.Cha}) the specific random choice procedure $\mathcal P$ used, as well as the main steps of the proof. The
orthogonal basis $e_1,\ldots,e_n$ of $W$ is formed as follows: $e_1$ is a random unit vector chosen uniformly in ${\mathbb R}^d$. The
second vector $e_2$ is a random unit vector from the space orthogonal to $e_1$, etc. For a fixed $x\in {\mathbb R}^d$ with $|x|_2=1$
the coordinates of the projection $\pi(x)\in W$ according to the basis $e_1,\ldots,e_n$ are independent random variables with the
expectation $\sqrt {1\over d}$. Application of the central limit theorem to $|\pi_W(x)|_2^2$, which is the sum of squares of $n$ coordinates,
provides, assuming $n\ge {20\over {1-p}} \ln m$, the following bound:
$P \{{1\over 2}|x|_2 \le \sqrt {d\over n} \ |\pi_W(x)|_2 \le 2|x|_2\} \ge 1-{{1-p}\over {m^2}}$. Finally, applying a union bound over
all $(_2^m)$ pairs of distances in $Z$, we obtain the required estimate.

\section{Entropy meets Dimensionality Reduction} \label{Entr.JLL}

In this section we prove Theorem \ref{main.short}. Let $Q$ be a relatively compact subset in a Hilbert space $V$ satisfying Assumption A, with the parameters $\beta, C$, and let $\e>0$ and $0<p<1$ be given. Put $\e_1={\e\over 6}$, and let $m=M+1$, where $M=M({\e_1},Q)$ is the covering number of $Q$ with $\e_1$-balls.

Let $x\in Q$ be an unknown vector. We have to show that for $n={20\over {1-p}} \ln m$, with probability at least $p$, the vector $x$ can be reconstructed with an error at most $\e$ from random $n$ measurements. The reconstruction algorithm proceeds in the following steps:

\smallskip

\noindent{\bf I. Preprocessing.}

\smallskip

1. Put $d=\left({RC\over \e_1}\right)^{1\over \beta}$, where $R>0$ is the radius of the minimal ball centered at $0\in V$, containing $Q$.

\smallskip

2. Take the $d$-dimensional subspace $V_d=span\{v_1,\ldots,v_d\}\subset V$, which exists according to the Assumption A, and let $\pi_{V_d}$ be the projection on
$V_d$. For any $z\in Q$ we have
\begin{align}\label{eq1}
|z-\pi_{V_d}(z)|_2 \le Cd^{-\beta}|z|_2 \le \e_1,
\end{align}
where the second inequality holds due to our choice of $d$, and since $|z|_2\le R$.

\smallskip

3. Pick at random a subspace $W \subset V_d$ of dimension $n$, and denote by $e_1,\ldots,e_n$ an orthonormal basis of it, and let $\pi_W:V\to W$ be the orthogonal projection of $V$ on $W$. Put $\hat \pi=\sqrt {\frac d n} \pi_W$, and fix the set of measurements, associated with $W$, as the set of linear functionals
$$
\lambda_i(z)=\langle\hat \pi(z),e_i\rangle ~, \quad i=1,\ldots,n ~,
$$
for each $z\in V$.

\smallskip

4. Let $M(\e_1,Q)$ be a minimal $\e_1$-cover of $Q$. Denote by $\{x_j\}_{j=1}^M \subset Q$ the centers of the covering $\e_1$-balls. Let $\tilde x_j = \pi_{V_d}(x_j), \ j=1,\ldots,M$, be the projections of $x_j$ to $V_d$. Notice that $\hat \pi (\tilde x_j)=\hat \pi (x_j) \in W$. Indeed,
$$
\hat \pi (\tilde x_j)=\sqrt {\frac d n} \pi_W\circ \pi_{V_d}(x_j)=\sqrt {\frac d n} \pi_W(x_j) = \hat \pi (x_j),
$$
since a composition of orthogonal projections on the
larger and on the smaller subspaces coincides with the projection on the smaller one. Denote the points
$$
y_j = \hat \pi (\tilde x_j)=\hat \pi (x_j)\in W ~,\quad j=1,\ldots,M .
$$

\smallskip

\noindent{\bf II. Reconstruction of $x$.}

\smallskip

For the unknown vector $x\in Q$ perform the measurements $\lambda_i(x), \ i=1,\ldots,n.$ Let $y=\hat \pi (x)\in W$. The coordinates
of $y$ in $W$ are exactly the measurements $\lambda_i(x)$, so $y$ is known. Take one of the points $y_j, \ j=1,\ldots,M$, say $y_{j_1}$ in the ball in
$W$ of radius $2\e_1$ around $y$. As we shall prove below, such points exist. We take $x_{j_1}\in Q$ as the desired $\e$-approximation of $x$.

\smallskip

\noindent{\bf Completion of the proof of Theorem \ref{main.short}.} It remains to show the with probability at least $p$, two assertions hold:

\smallskip

a. The existence of points $y_j$ in the ball in $W$ of radius $2\e_1$ around $y$.

\smallskip

b. For each one of them, say, for $y_{j_1}$, we have $|x-x_{j_1}|_2\le \e.$

\medskip

Consider the point $\tilde x=\pi_{V_d}(x)\in V_d$. We have $y=\hat \pi (x)= \hat \pi (\tilde x) \in W$, by the same reasoning as above.
Now, define $Z\subset V_d$ as
$$
Z=\{\tilde x, \tilde x_1 \ldots, \tilde x_M\}.
$$
So $Z\subset V_d$ is a finite set of cardinality $m=|Z|=M+1$. Let us recall that $\tilde x$, as well as $x$, is unknown, while
$\tilde x_j$'s are known points. Now by Theorem \ref{JLLT} with probability $p$ for any $u,v \in Z$ we have
\begin{align}\label{eq3}
\frac 1 2 |u-v|_2 \le \sqrt {\frac d n } \Big|\pi_W(u)-\pi_W(v)\Big|_2 \le 2|u-v|_2.
\end{align}

Notice that there exists $j_0$ such that
\begin{align} \label{eq2}
|\tilde x - \tilde x_{j_0}|_2 \le \e_1.
\end{align}
Indeed, $\|\pi_{V_d}\|=1$, and $x\in Q$ belongs to one of the covering balls, so the center of this ball $x_{j_0}$ satisfies 
$|x-x_{j_0}|_2 \le \e_1$. Thus, by the definition of $y,y_{j_0}$, the upper bound of \eqref{eq3} and \eqref{eq2} we have
$$
|y-y_{j_0}|_2 \le 2\e_1.
$$
This proves the assertion (a) above, i.e. the existence of points $y_j$ in the ball in $W$ of radius $2\e_1$ around $y$.

Let $y_{j_1}$ be one of such points. It remains to show that $|x-x_{j_1}|_2 \le \e$. Indeed,
$$
|x-x_{j_1}|_2 \le |x-\tilde x|+|\tilde x-\tilde x_{j_1}|+|\tilde x_{j_1}-x_{j_1}| .
$$

The first and the third terms, by \eqref{eq1}, do not exceed $\e_1$. Now, since $y_{j_1}$ is in the ball in $W$ of radius $2\e_1$
around $y$, we have $|y-y_{j_1}|_2\le 2\e_1.$ Applying the lower bound of \eqref{eq3}, we conclude that
$|\tilde x-\tilde x_{j_1}|_2 \le 4\e_1$. Finally,
$$
|x-\tilde x_{j_1}|_2 \le 6\e_1 = \e,
$$
which completes the proof of Theorem \ref{main.short}. $\square$

\smallskip

\begin{rem} \label{fin.acc}
Let us assume that our measurements $\lambda_v (x)=\langle x,v \rangle$ are normalized, i.e. $|v|_2=1$, and have a finite accuracy $\delta >0$. In order to keep the reconstruction accuracy $\e$ we have (according to our proof of Theorem \ref{main.short} given above) to assume that $\delta \le \frac{\e}{\sqrt d}$, where $d$ is the dimension of the ambient space. This is because the linear mapping $\hat \pi$ providing the dimensionality reduction in the Johnson-Lindenstrauss lemma, Theorem \ref{JLLT} above, is an orthogonal projection amplified by, essentially, $\sqrt d$. Hence this is the expected amplification factor of the measurement noise. The dimension $d=({1\over \e})^{1\over \beta}$ in typical examples is of order $\frac 1 \e$.  Altogether we have to require $\delta \asymp \e^{3\over 2}$, where `` $\asymp$ '' denotes equivalence up to constants, not depending on $\e$. Whether such a high accuracy of the samples is indeed necessary, is one of the open problems in our approach. Another open problem concerns the true necessity of the random choice of the measurements, which seems to be unavoidable in our specific algorithm above.
\end{rem}

\smallskip

\begin{rem}
We do not provide in this paper lower bounds for sampling entropy. Simple information considerations show that
$N(\e,Q)$ cannot be much smaller than $H(\e,Q)$. Indeed, let us estimate the maximal amount of information provided by one measurement with
the range $[-1,1]$ and a possible error of order $\delta$. With this accuracy we can surely identify only the measurement values quantified up
to $2\delta$, and this gives us $\log ({1\over \delta})$ binary bits of information. Thus, $N$ measurements provide at most $N \log ({1\over \delta})$
bits. But to specify any given $x\in Q$ up to accuracy $\e$ we need exactly $H(\e,Q)$ bits, by the definition of the Kolmogorov entropy, no matter
what is the source of the information. Therefore, independently of the specific reconstruction method we apply, we must have
$N \log ({1\over \delta})\ge H(\e,Q)$, or

\be\label{low.bd}
N(\e,Q) \ge H(\e,Q)\log^{-1} ({1\over \delta}).
\ee
Under the assumption $\delta \asymp \e^{3\over 2}$, as in Remark \ref{fin.acc}, we get $N(\e,Q) \ge {2\over 3}H(\e,Q)\log^{-1} ({1\over \e}).$
\end{rem}

\smallskip

\begin{rem}
For exact measurements, in a fixed finite dimensional setting, accurate bounds for sampling entropy are provided
in \cite{Ver1} via the ``mean width'' of $Q$. It would be interesting to compare the results of the present paper and of \cite{Ver1} in a
situation where the dimension of the ambient space grows as the required accuracy bound $\e$ tends to zero.
\end{rem}

\smallskip

\section{Main Examples}\label{EXa}

In this section we state some known and prove some new results on the Kolmogorov entropy of important functional classes. Through Theorem
\ref{main.short} we obtain estimates for the corresponding sampling entropy, comparing them with the known values, when available.

\subsection{Classical regularity spaces}

The following result has been established in \cite{Kol1,Kol2,Kol.Tih,Zak}:

\bt
Let $C^k(n,K)\subset C_0(I^n)$ be a (relatively compact) subset of $C_0(I^n)$ consisting
of all $C^k$-smooth functions on $I^n$ with all the derivatives up to order $k$ bounded by
$K$. Then
$$
H(\e,C^k(n,K))\asymp ({1\over \e})^{n\over k} .
$$
Moreover, let $A=A(\eta,K) \subset C_0([-\pi,\pi])$ consist of periodic real analytic functions on
$[-\pi,\pi]$, extendable into a complex $\eta$-neighborhood of $[-\pi,\pi]$ and
uniformly bounded there by $K$. Then
$$
H(\e,A)\asymp \log^2({1\over \e}) .
$$
\et
Kolmogorov's results have been obtained for the $C_0$ ambient metric. However, for the classes
$C^d(n,K)$ and $A(\eta,K)$ consisting of functions with uniformly bounded derivatives, the
$L^2$ and $C_0$ norms are equivalent, so the results remain true also in $L^2$ ambient norm.

\smallskip

Now the following corollary of our general Theorem \ref{main.short} recovers, in principle, the
classical sampling accuracy for smooth functions. For analytic ones it gives less sharp result: Indeed,
sampling with the first $N\asymp \log({1\over \e})$ Fourier coefficients reconstructs analytic function
with an accuracy $\e$.
\bc
$$
N(\e,p,C^k(n,K))\le {{C(n,k)}\over {1-p}}\biggl({4\over \e}\biggr)^{n\over k} ~, \quad N(\e,p,A(\eta,K))\le {{C_1(\eta,K)}\over {1-p}}\log^2\biggl({1\over \e}\biggr).
$$
\ec

\subsection{Piecewise Smooth Functions}\label{P.Sm.Fns}

In this section we show that the metric entropy of piecewise smooth functions is essentially
the same as for their smooth counterparts. This is intuitively clear in one variable: In
addition to the ``smooth information'' we have just to specify the positions of the jumps,
which requires a negligible amount of additional information. In several variables we have
to be more careful: Complexity of the boundaries of the smooth regions can significantly
contribute to the overall entropy.

In what follows we shall consider piecewise $C^k$-smooth functions $f(x_1,\dots,x_n)$ defined
on the torus $T^n$ which we identify with the $n$-dimensional box $I^n=[-\pi,\pi]^n$ with the
opposite faces glued together.

\smallskip

\noindent {\bf Assumptions B}.

\smallskip

We shall assume that the jumps of $f$ occur at certain $C^k$-smooth,
pairwise disjoint, compact hypersurfaces $\Sigma_j \subset T^n, \ j=1,\dots,K.$ The
complement in $T^n$ of $\Sigma=\cup_{j=1}^K \Sigma_j$ is a union of certain $n$-dimensional
manifolds $D_s$ whose boundaries $\partial D_s$ are unions of some of $\Sigma_j$. So we
assume that $f=f_s$ on each $D_s$, with $f_s$ being $C^k$-smooth functions on $D_s$.

The following parameters of $f$ will be essential in formulation of our result: First of all,
let $\rho_1$ denote the minimal distance between the jump submanifolds $\Sigma_j$. Next,
denote for each point $x\in \Sigma$ by $n(x)$ the unit normal vector to $\Sigma$ at $x$.
Consider a new orthogonal coordinate system $(y_1,\dots,y_n)$ at $x$ such that the axis $Oy_n$
is parallel to $n(x)$. Locally the hypersurface $\Sigma$ near $x$ can be represented by
the equation $y_n=\eta(y_1,\dots,y_{n-1})$, with $\eta$ a $C^k$-smooth function. We shall
assume that for each point $x\in \Sigma$ this representation is valid in a ball $B$ of a
radius at least $\rho_2$ in the coordinates $y_1,\dots,y_{n-1}$, with all the derivatives
of $\eta$ up to order $k$ uniformly bounded by $K_1$ in $B$.

The constants $\rho_2$ and $K_1$ are closely related to the injectivity radius of $\Sigma$,
and they can be estimated through the upper bound on its curvature.

Next, we assume that all the derivatives of the smooth pieces $f_s$ of $f$ up to order $k$
are uniformly bounded by $K_2$ in the corresponding domains $D_s$. This assumption implies
the upper bound $2K_2$ for the jumps of all the derivatives of $f$ along $\Sigma$.

\bd\label{pCk}
The class $PC^k(n)=PC^k(n,\rho_1,\rho_2,K_1,K_2)\subset L^2(T^n)$ consists of all the
piecewise $C^k$-smooth functions on the torus $T^n$, satisfying Assumptions B with the
parameters $\rho_1,\rho_2,K_1,K_2$.
\ed
\bt\label{pCke}
The Kolmogorov $\e$-entropy of $PC^k(n) \subset L^2(T^n)$ satisfies
$$
H(\e,PC^k(n))\asymp ({1\over \e})^m ,
$$
where $m={1\over k}$ for $n=1$, and $m={{2(n-1)}\over k}$ for $n\ge 2$.
\et
\pr
We shall construct explicitly an $\e$-net in the set $PC^k(n)\subset L^2(T^n)$ containing
approximately minimal number of elements. Let us start with dimension $n=1$. Assumptions B
imply the uniform bound $s$ on the number of the jump points of $f\in PC^k(1)$.
First we quantize, up to accuracy $\e^2$, the positions of these jump points. $\e^2$ appears
here, instead of $\e$, since a shift of the jump points to a distance of order $h$ produces
$L^2$-norm of order $\sqrt h$ for the difference of the corresponding piecewise smooth functions.

There are at most $({{2\pi}\over \e})^{2s}$ quantized jump configurations. For each configuration
we construct an $\e$-net of $C^k$ functions on each continuity segment $\Delta_j$. Taking all the
jump configurations, and picking independently elements of $\e$-nets on each continuity segment
$\Delta_j$, we construct an $\e$-net in $L^2([-\pi,\pi])$ for the entire set $C^k(1)$. The number
of the elements in this net is the product of the number of quantized jump configurations, and of
the number of the elements in the nets for each of the continuity segments. After taking
logarithm we get, by the Kolmogorov result, that $H(\e,PC^k(1))\asymp ({1\over \e})^{1\over k}.$
Indeed, the contribution of the logarithm of the number of quantized jump configurations is
at most $2s \log ({{2\pi}\over \e})$, which is negligible with respect to the main asymptotic
rate.

In several variables, we repeat the same construction. The jump manifolds $\Sigma$ are, by
Assumptions B, $C^k$ smooth $(n-1)$-dimensional submanifolds in $T^n$, with uniformly bounded
derivatives up to order $k$ (in appropriate coordinate charts $\eta$). By Kolmogorov's result
in dimension $n-1$, after quantization with the step $\e^2$ and taking logarithm, we get
$\log M_0 \asymp ({1\over \e})^{{{2(n-1)}\over k}}$ for the number $M_0$ of the quantized
jump configurations. For the number $M_j$ of the elements in the $\e$-net on each continuity
domain $D_j$ we have $\log M_j \asymp ({1\over \e})^{{{n}\over k}}.$ But now, in contrast with
the one-dimensional case, we have ${{2(n-1)}\over k}>{{n}\over k}$ for $n\ge 3$. So the
asymptotic behavior is determined by $m=\max ({{2(n-1)}\over k},{{n}\over k})$, which is equal
to ${n\over k}$ for $n=1$, and to ${{2(n-1)}\over k}$ for $n\ge 2$. This completes the
proof of Theorem \ref{pCke} $\square$

\smallskip

By Theorem \ref{main.short} we conclude that for piecewise smooth functions in one and two variables
the sampling entropy is the same as for their smooth counterparts, while starting with three
variables it becomes strictly larger. This is an artifact of the $L^2$ norm: It requires double
accuracy in positioning of the jumps. In an appropriate pseudo-norm, combining $C^0$-distance of smooth pieces, 
and $C^0$-distance of the jump submanifolds, the Kolmogorov $\e$-entropy, and hence the sampling
$\e$-entropy, remain of order $({1\over \e})^{{{n}\over k}}.$

This fact was conjectured in various forms in many publications in the field. Let us consider in somewhat more detail one specific
setting of this conjecture.

Periodic $C^k$ functions on $[-\pi,\pi]$ are approximated by their partial Fourier sums up to the frequency $N$ with the accuracy
of order $N^{-k-1}$. For real analytic functions the approximation accuracy is of order $q^{-N}, \ q<1.$ From this fact it follows that
the deterministic sampling $\e$-entropy of these classes is of order $({1\over \e})^{1\over {k+1}}$ and $\log ({1\over \e})$, respectively.
However, if discontinuity points appear, this rate drops to $1\over N$, no matter what is the signal regularity between the jumps. It can
be shown that if the position of the jumps is free and not known in advance, any linear method cannot provide better accuracy than
$1\over N$ (see, for example, \cite{Eti.Sar.Yom}).

Fourier reconstruction of piecewise smooth functions is important in many theoretical and applied problems, from numerical
solution of non-linear PDE's developing shock waves, to medical imaging. Many recent publications have been devoted to various
methods of such reconstruction (see \cite{Adc.Han,Bat.Yom1,Eck,Gel.Tad,Hry.Gro,Kve} and references therein for a very
partial sample). While it is known that the ``smooth'' sampling rate can be achieved inside the continuity intervals, the problem
of accurate reconstruction of the positions of the jump points (and the signal nearby) turns out to be more difficult. See, however,
\cite{Lip.Lev} where piecewise smooth functions which are continuous, but have discontinuous derivatives, are reconstructed from point-wise
samples with full accuracy). 

It was conjectured by Eckhoff and other authors (\cite{Eck,Kve} and references therein) that the ``smooth'' accuracy rate $N^{-k-1}$
of Fourier reconstruction of univariate $C^k$ functions from their first $N$ Fourier coefficients can be achieved also in a
piecewise smooth case (including the positions of the jumps), if we allow non-linear processing of the Fourier data. In spite of many
attempts, this conjecture remained largely open. In \cite{Bat.Yom1} we have shown that a modification of the original Eckhoff's method
(based on algebraic reconstruction of the positions and amplitudes of the jumps) produces ``half of the accuracy'', $N^{-{k\over 2}}$.
Recently Eckhoff's conjecture was finally settled in \cite{Bat}, by a further significant improvement of the original Eckhoff's
method. For the the deterministic sampling $\e$-entropy of piecewise $C^k$-smooth functions this result gives the ``smooth order''
$({1\over \e})^{1\over {k+1}}$. However, the analytic version of this problem remains open. So, the result of the next section is
apparently, a first confirmation of the ``Eckhoff-type'' conjecture for univariate piecewise-analytic functions.

\subsection{Piecewise Analytic functions}

To simplify the presentation, we restrict ourselves only to the case of one variable. Similar result
remains true in a multidimensional situation.

\bt Let $PA=PA(\kappa,\eta,K)\in L^2([-\pi,\pi])$ be the set of piecewise real analytic
functions on $[-\pi,\pi]$ with at most $\kappa$ discontinuity points, such that each
their analytic piece can be extended into a complex $\eta$-neighborhood of the
corresponding real segment, and is bounded there by $K$.
Then
$$
H(\e,PA)\asymp \log^2({1\over \e}) .
$$
\et
\pr
Exactly as for Theorem \ref{pCke}, one-dimensional case. $\square$

\smallskip

By Theorem \ref{main.short} we conclude that for piecewise analytic functions in one variable
the sampling entropy is of order $\log^2({1\over \e}).$ It is an open question whether this
bound can be improved to the classical order of $\log({1\over \e}).$

\subsection{Finite-Dimensional Warping}

Let a finite-dimensional family of transformations
$\Psi_\tau:[-\pi,\pi]\rightarrow [-\pi,\pi]$ be given, with the parameter
$\tau=(\tau_1,\dots,\tau_s)\in I^s= [0,1]^s$. We shall assume that $\Psi_\tau(x)$ is
Lipschitzian, both in $x$ and in the parameter $\tau$.

\smallskip

Let a relatively compact $Q\in L^2([-\pi,\pi])$ be given, with all $f\in Q$ being
Lipschitzian, with the same Lipschitz constant. Define
$Q\circ \Psi:=\{g=f\circ \Psi_\tau, \ f\in Q, \tau \in I^s \}$. Since $\Psi_\tau$ are
Lipschitzian, $Q\circ \Psi$ is also a a relatively compact subset in $L^2([-\pi,\pi])$.

\bt\label{Warp}
For each $\e>0$ we have $H(\e,Q\circ \Psi) \le H(\e,Q)+C\log({1\over \e})$. If
$\Psi_\tau$ are homeomorphisms of $[-\pi,\pi]$ then $H(\e,Q\circ \Psi) \le H(\e,Q)$.
\et
\pr
This is a very special case of the results proved in \cite{Kol1,Kol.Tih,Vit} on the metric
entropy of composition classes. Notice that in general $H(\e,Q\circ \Psi)$ may be strictly
smaller than $H(\e,Q)$. This happens, in particular, if, $\Psi$ maps the entire interval
$[-\pi,\pi]$ into one point. $\square$

\smallskip

In particular, let $C^k,A$ denote the classes of $C^k$-smooth and real analytic functions on $[-\pi,\pi]$

\bc
For $\Psi_\tau$ a finite-dimensional family of homeomorphisms we have
$$
H(\e, C^k \circ \Psi)\asymp ({1\over \e})^{1\over k} ~, \quad H(\e, A \circ \Psi)\asymp \log^2 ({1\over \e}) .
$$
\ec

By Theorems \ref{Warp} and \ref{main.short} we conclude that the sampling entropy of $C^k \circ \Psi, \ A \circ \Psi$
remains the same as for their regular counterparts, while the regularity may be completely lost.

\subsection{Adding Finite-Dimensional Linear Combinations}

It is easy to show that adding to the functional classes as above linear combinations of a finite number of fixed
Lipschitzian functions, we do not change the order of their Kolmogorov entropy. By Theorem \ref{main.short} the sampling
complexity remains the same as for their regular counterparts, while the regularity may be completely lost.

\section{Some questions and comments}

1. Is it possible to take into account a ``regular'' geometry of the $\e$-nets that appear
in the functional classes discussed in this paper, in particular, for analytic functions? It
is well known that the dimension bound in Johnson-Lindenstrauss lemma is sharp, up to a factor
$\log ({1\over \sigma})$, where $\sigma$ is the almost-isometry distortion allowed (See \cite{Alo}).
In our applications $\sigma={1\over 2}$ is fixed, but we encounter a gap of order $\log ({1\over \e})$
in the projection dimension for our $\e$-nets, in particular, in the case of analytic functions.
Are these two effects related? Compare also bound \ref{low.bd} above for finite accuracy measurements.

\smallskip

2. What is the true accuracy rate in sampling of piecewise-analytic functions?
It is known that inside the continuity intervals the classical reconstruction rates from
Fourier coefficients can be recovered (\cite{Gel.Tad}).

For piecewise smooth functions the same gap in the upper and lower bounds appears as in the piecewise analytic case.
However, it does not change the asymptotic recovery rate, because of a much higher complexity of piecewise smooth
functions.

\smallskip

3. What is the ``true'' order of the allowed measurements errors in the results above? As it was mentioned
in Remark \ref{fin.acc}, we can modify the proofs above to include measurements of a finite accuracy, but the
order of the allowed error is $\e^{3\over 2}$, which may be too high.

\vskip1cm

\bibliographystyle{amsplain}

\end{document}